\documentclass{amsart}
\usepackage{amsmath,amsthm,amssymb,color}
\usepackage{hyperref}
\newtheorem{proposition}[equation]{Proposition}

\newtheorem{lemma}[equation]{Lemma}
\newtheorem{corollary}[equation]{Corollary}
\newtheorem{definition}[equation]{Definition}

\newtheorem{assumption}[equation]{Assumption}
\theoremstyle{remark}
\newtheorem{example}[equation]{Example}

\newcommand\inv{^{-1}}
\newcommand\lexp[2]{\kern\scriptspace\vphantom{#2}^{#1}\kern-\scriptspace#2}

\DeclareMathOperator\rightlcm{{\mathrm right-lcm}}

\DeclareMathOperator\leftgcd{{\mathrm left-gcd}}

\newcommand\Ps{{P^{(s)}}}
\author{Fran\c cois Digne and Jean Michel}
\address[F.~Digne]{Laboratoire Ami\'enois de Math\'ematique Fondamentale et
Appliqu\'ee, CNRS UMR 7352, Universit\'e de Picardie-Jules Verne,
80039 Amiens Cedex France.}
\email{digne@u-picardie.fr}
\urladdr{www.lamfa.u-picardie.fr/digne}
\address[J.~Michel]{Institut Math\'ematique de Jussieu -- Paris rive
gauche, CNRS UMR 7586, Universit\'e de Paris, B\^atiment Sophie Germain,
75013, Paris France.}
\email{jean.michel@imj-prg.fr}
\urladdr{webusers.imj-prg.fr/$\sim$jean.michel}


\title{Ribbons in Garside monoids} 
\begin{document}
\begin{abstract}
We expound the properties of ribbons in a setting 
which  is  general  enough  to encompass spherical Artin
monoids and dual braid monoids of well-generated complex reflection groups.
We generalize to our
setting results on parabolic subgroups of spherical Artin group of Godelle
\cite{Godelle},  Cumplido \cite{Cumplido}  and others.
\end{abstract}
\maketitle
This short note expounds the properties of ribbons in a setting (Assumptions
\ref{assumption} and \ref{pow delta central})  which  is  general  enough  to encompass spherical Artin
monoids and dual braid monoids of well-generated complex reflection groups.
We show in particular that this approach can recover, and generalize to our
setting, results on parabolic subgroups of spherical Artin group of Godelle
\cite{Godelle},  Cumplido \cite{Cumplido}  and others.  We have strived for
the  text to be  self-contained apart from  basic results on Garside theory
for which we refer to \cite{livre}. The assumptions under which we study
ribbons are slightly stronger than those in  \cite[VIII \S 1.4]{livre} but
allow us to develop them  much further.

Let   $M$  be  Garside  monoid;  that  is, $M$  is  Noetherian,  left-  and
right-cancellative,  any pair  of elements  have left-  and right-gcd's and
lcm's,  and there exists a Garside element $\Delta$, an element whose left-
and right-divisors coincide and generate $M$ (see \cite[I 2.1]{livre}).

We  recall  that  a  factor  is  any  left-divisor  of  a right-divisor (or
equivalently  a right-divisor of a  left-divisor) and that Noetherian means
there  is no infinite  sequence where each  term is a  proper factor of the
previous one.

\begin{definition}
A  standard parabolic submonoid $P$ of $M$ is a submonoid closed by factors
such that any element of $M$ has a maximal left-divisor (the left $P$-head)
in $P$ and a maximal right-divisor (the right $P$-head) in $P$.
\end{definition}
This is the particular case, for a monoid, of Definition \cite[1.30 VII]{livre}
which is in the context of a Garside category.

We denote by $H_{P}(g)$ the left $P$-head of $g\in M$ and
we define $T_{P}(g)$ (the "$P$-tail") by $g=H_{P}(g)T_{P}(g)$. We say that $g$
is $P$-reduced if $H_P(g)=1$ or equivalently $g=T_P(g)$.

Any  intersection of standard parabolic  submonoids is a standard parabolic
submonoid  \cite[VII 1.35]{livre}. It is  clear that any standard parabolic
submonoid  of  a  standard  parabolic  submonoid  is  a  standard parabolic
submonoid.

Since  $M$  is  Noetherian,  it  is  generated  by the set $S$ of its atoms
(elements  which  have  no  proper  nontrivial  divisors)  and  a  standard
parabolic submonoid being closed by factors is also generated by its atoms.

\begin{lemma}\label{garside element}
A standard parabolic submonoid $P$ has a Garside element $\Delta_P$ given by $\Delta_P=
H_{P}(\Delta)$.
\end{lemma}
\begin{proof}  First $P$  is generated  by the  left-divisors of $\Delta_P$
since  it is  generated by  its atoms  which are  left-divisors of $\Delta$
hence  of $\Delta_P$. Now the right- and left- divisors of $\Delta$ are the
same, thus $\Delta\succcurlyeq\Delta_P$. Hence
$H^\succcurlyeq_P(\Delta)\succcurlyeq\Delta_P$,   if  $H^\succcurlyeq_P(g)$
denotes  the  right  $P$-head  of  $g$.  Exchanging  right  and left, that is
$H_P$ and $H^\succcurlyeq_P$, we get
symmetrically     $H^\succcurlyeq_P(\Delta)\preccurlyeq\Delta_P$.     Hence
$\Delta_P$  is  a  factor  of  itself.  It  cannot  be  a  strict factor by
Noetherianity, thus
$\Delta_P=H_P^\succcurlyeq(\Delta)$.  The  left-divisors  of $\Delta_P$ are
left, thus right-divisors  of $\Delta$, and since they are in $P$  they are
right-divisors of $H_P^\succcurlyeq(\Delta)=\Delta_P$. Conversely
right-divisors of $\Delta_P$ are right, thus left-divisors of $\Delta$ and
being in $P$ are left-divisors of $\Delta_P$.
\end{proof}

\begin{assumption}\label{assumption}
We assume that for any subset $S_1\subset S$, the right-lcm of $S_1$ is the 
Garside element given by Lemma \ref{garside element} of the
smallest standard parabolic submonoid containing $S_1$.
\end{assumption}
In the following we assume that Assumption \ref{assumption} holds. It clearly
holds in a spherical Artin monoid.

\begin{proposition}The dual braid monoid $M$ of a well generated complex
reflection group $W$ satisfies Assumption \ref{assumption}.
\end{proposition}
\begin{proof}
From the fact that any decomposition of $\Delta$ is obtained from one of them
by the Hurwitz action (see \cite[Proposition 7.6 and Proposition 8.5]{Bessis}), it follows that the same holds for any
simple.  It follows from that that any factor of a simple $w$ is a left-divisor of $w$.
It follows that any simple $w$ is a common right multiple of the atoms which left-divide
it. If the right-lcm of these atoms was a strict divisor of $w$,
there would be an atom $s$ whose square is a factor of $w$, thus of $\Delta$.
This does not exist, see \cite[Property (iii) of Proposition 2.1 (M2)]{DMM}.

Thus any simple $w$ is the left-lcm of the set $S_w$ of atoms which left divide it. 
The submonoid $P_w$ generated by these atoms is standard parabolic since it is
clearly stable by factor, and every element has a $P$-head : to see this it is
sufficient to check it for a simple $v$ by \cite[VII, 1.25]{livre}, in which
case the left-gcd of $v$ and $w$ is the head.

Finally, for any set $S_1$ of atoms, the right-lcm $w$ of this set defines
thus a parabolic $P_w$ which is clearly the minimal parabolic subgroup
containing $S_1$.
\end{proof}
If $b,g\in M$ and $g\preccurlyeq bg$ we say that the conjugate $b^g$,
equal to $g\inv b g$, is defined. For a set $P\subset M$ we say that $P^g$
is defined if $p^g$ is defined for any $p\in P$.
\begin{definition} For $P$ a standard parabolic submonoid
we call $P$-ribbon a $P$-reduced element $g\in M$
such that $P^g$ is defined. We denote such a ribbon by $P\xrightarrow{g}$.
\end{definition}
Let $P$ be a standard parabolic submonoid  with Garside element $\Delta_P$ and
$s$ be an atom of $M$ which is not in $P$;
we define $\Ps$ as the smallest
standard parabolic submonoid  containing $P$ and $s$.
Since $\Delta_\Ps$ is left-divisible by
$\Delta_P$ the following makes sense:
\begin{definition}\label{v(s,M)} For $P$ a standard parabolic submonoid
and $s\in S, s\notin P$, we define $v_{s,P}$ by $\Delta_\Ps=\Delta_P v_{s,P}$.
\end{definition}

Recall that conjugation by $\Delta$ is an automorphism of $M$
(see for example \cite[V 2.17]{livre}).

\begin{lemma}
For $P$ a standard parabolic submonoid and for $s\in S$, $s\notin P$ the element
$v_{s,P}$ is a $P$-ribbon. 
\end{lemma}
\begin{proof}
 We have $\Delta_P=H_P(\Delta_\Ps)$ since any divisor in $P$ of
 $\Delta_\Ps$ divides $\Delta$
hence divides $\Delta_P$. Hence $v_{s,P}$ is $P$-reduced. And
 $P^{v_{s,P}}$ is defined since the conjugation by
$v_{s,P}$ is the composition of the inverse of the conjugation by $\Delta_P$, which is an
automorphism of $P$, and the conjugation by $\Delta_\Ps$, which is an
automorphism of $\Ps$.
\end{proof}
\begin{lemma}\label{atoms of the ribbon category}
If $g$ is a $P$-ribbon and $s\in S$ left-divides $g$ then
$v_{s, P}\preccurlyeq g$.
\end{lemma}
\begin{proof}
We have  $\Delta_P g=g\Delta_P^g$.  
Thus  $s\preccurlyeq \Delta_Pg$, thus 
 $\Delta_\Ps\preccurlyeq\Delta_P  g$ by Assumption \ref{assumption},  thus
 $v_{s,P}=\Delta_P\inv\Delta_\Ps\preccurlyeq g$.
\end{proof}
Note that the converse of this lemma, that is $s\preccurlyeq v_{s,P}$, holds in
an Artin monoid but not necessarily in a dual braid monoid, see
Example \ref{dual monoid}.
\begin{proposition}\label{category}Let $P$ be a standard parabolic submonoid; then
\begin{itemize}
\item  If $P\xrightarrow{g}$ is a ribbon then $P^g$ is a standard parabolic
submonoid,  and $x\mapsto  x^g$ is  a monoid isomorphism $P\xrightarrow\sim
P^g$.
\item A $P$-ribbon is a product of elements $v_{s,P'}$ where $P'$ are standard
parabolic submonoids conjugate to $P$.
\item If $P\xrightarrow{g}$ is a ribbon then the conjugate by $g$ of the atoms
of $P$ are the atoms of $P^g$.
\item If $P\xrightarrow{g}$ is a ribbon then $\Delta_P^g=\Delta_{P^g}$.
\end{itemize}
\end{proposition}
\begin{proof}
Let  $s$  be  an  atom  dividing  $g$,  then  $v_{s,P}$ conjugates $P$ on a
 standard  parabolic submonoid since the  conjugation by $\Delta_\Ps$ is
an  automorphism of $\Ps$, hence conjugates $P$ on a standard parabolic
submonoid  of $\Ps$ thus of $M$, and the conjugation  by $\Delta_P$ conjugates  $P$ on itself. If
$g_1=  v_{s,P}\inv g$, it follows  that $P^{v_{s,P}}\xrightarrow{g_1}$ is a
ribbon ($g_1$ is $P^{v_{s,P}}$-reduced since for an atom $s$ if
$s^{v_{s,P}}\preccurlyeq g_1$ then $s\preccurlyeq g$). By Noetherian induction we get the first two items.

The third item is an immediate consequence of the first one, and the fourth
also by Assumption \ref{assumption}.
\end{proof}
In view of the third item above, the definition of a $P$-ribbon could 
instead of asking that $P^g$ is defined just ask that the conjugate of any 
atom of $P$ by $g$ is defined.

\begin{lemma}\label{lcm(t,v)}
If $P$ is a standard parabolic submonoid, if $s$ is an atom of $P$ and
$P\xrightarrow g$ is a  ribbon,
then the right-lcm of $s$ and $g$ is $sg$.
\end{lemma}
\begin{proof}
Since $g$ is a $P$-ribbon it conjugates $s$ into some atom $s'$. The
right-lcm of $s$ and $g$ left-divides
$s g=gs'$ and is a strict right-multiple of $g$ (since $g$ is $P$-reduced), 
hence is equal to $gs'$.
\end{proof}

\begin{lemma} \label{h<TP(g)} Let $P\xrightarrow g$ be a
ribbon;  for $x\in P$, $y\in M$, it is equivalent that  $g\preccurlyeq xy$ or that
$g\preccurlyeq y$.
\end{lemma}
\begin{proof}
If $y=gy'$, then $xy=gx^gy'$ so that $g\preccurlyeq xy$. To prove the converse, by
Noetherian induction on $x$, it is sufficient to prove that if  an atom $s\in P$ left-divides
$gh$ for some $h$, then $s^g\preccurlyeq h$.
By Lemma \ref{lcm(t,v)}, the right-lcm of $s$ and $g$ is $sg=gs^g$.
Thus $s\preccurlyeq gh$ is equivalent to
$gs^g\preccurlyeq gh$ which is finally equivalent to $s^g\preccurlyeq h$.
%
\end{proof}
\begin{lemma}\label{I-head preserved}
Let $P\xrightarrow g$ be a ribbon and
let $h\in M$. Then $T_{P}(gh)=gT_{P^g}(h)$ and
$H_{P}(gh)^g=H_{P^g}(h)$. In particular if $g$ is a $P$-ribbon and $h$ is a $P^g$-ribbon,
then $gh$ is a $P$-ribbon.
\end{lemma}
\begin{proof} Let $s$ be an atom of $P$ and set $s':=s^g$. Both formulae clearly
follow if we show that it is equivalent that $s\preccurlyeq gh$ or that
$s'\preccurlyeq h$. This is true by the proof of Lemma \ref{h<TP(g)}.
\end{proof}

It follows from Lemma \ref{I-head preserved} and the first item of Proposition
\ref{category} that the ribbons form a category.  The second item of Proposition
\ref{category} shows that the atoms of this category are the
$P\xrightarrow{v_{s,P}}$ which are not strict multiple of some other one.
This is always the case in a spherical Artin monoid but not in a dual braid
monoid as the following example shows.
\begin{example}\label{dual monoid}
We consider the Coxeter Group of type $A_4$ identified
with the symmetric group on 5 letters. 
We number the transpositions $s_1,\ldots,s_{10}$
in the order $(1,2), (2,3), (3,4), (4,5), (1,3), (2,4), (3,5), (1,4), (2,5), (1,5)$.
We choose as Coxeter element the product $s_1\cdots s_4$. The atoms of the
corresponding dual monoid are in one to one correspondence
with the transpositions.  We denote the atoms again by $s_1,\ldots s_{10}$.
Let $P$ be the standard parabolic subgroup
generated by the atom $s_5$; then $\rightlcm(s_5,s_6)=s_5s_1s_3$, thus 
$v(s_6,P)=s_1 s_3$ and $\rightlcm(s_5,s_1)=s_5s_1$ so that $v(s_1,P)=s_1$.

Note that $s_6$ does not divide $v(s_6,P)=s_1 s_3$.
\end{example}

\begin{lemma}\label{gcd lcm defined}
If $b,g,g'\in M$ and $b^g$ and $b^{g'}$ are defined, then so
are $b^{\leftgcd(g,g')}$ and $b^{\rightlcm(g,g')}$.
\end{lemma}
\begin{proof}
The conditions for $b^g$ and $b^{g'}$ to be defined are $g\preccurlyeq bg$
and $g'\preccurlyeq bg'$. It follows easily that
$\leftgcd(g,g')\preccurlyeq b\leftgcd(g,g')$ and
$\rightlcm(g,g')\preccurlyeq b\rightlcm(g,g')$.
\end{proof}
For $g\in M$ we denote by $H(g)$ the first term of its Garside normal form,
equal to $\leftgcd(g,\Delta)$ and define $T(g)$ by $g=H(g)T(g)$.
\begin{proposition}
If $P\xrightarrow g$ is a ribbon, so are all the terms
of its Garside normal form.
\end{proposition}
\begin{proof}
It is sufficient to prove that $P\xrightarrow{H(g)}$ and
$P^{H(g)}\xrightarrow{T(g)}$ are ribbons.

For the first fact, we have $H_{P}(H(g))=1$ since $H_{P}(g)=1$, and
$P^{H(g)}$ is defined by  Lemma \ref{gcd lcm defined}, since 
$H(g)=\leftgcd(g,\Delta)$.

For  the second fact,  since $(P^{H(g)})^{T(g)}=P^g$ is defined,
it is sufficient to show that $T(g)$ is $P^{H(g)}$-reduced. By
Lemma \ref{I-head preserved} with $H(g),T(g)$ for $g,h$, we get that
$H_{P^{H(g)}}(T(g))=H_P(H(g)T(g))^{H(g)}=1$.
\end{proof}
\begin{proposition}
If $P\xrightarrow g$ and $P\xrightarrow{g'}$ are in the ribbon category, then 
also the map $P\xrightarrow{\leftgcd(g,g')}$
\end{proposition}
\begin{proof}
It is clear that $H_{P}(\leftgcd(g,g'))=1$, and
$P^{\leftgcd(g,g')}$ is defined by lemma \ref{gcd lcm defined}.
\end{proof}
\begin{proposition}
If $P\xrightarrow g$ and $P\xrightarrow{g'}$ are in the ribbon category, then 
also the map $P\xrightarrow{\rightlcm(g,g')}$.
\end{proposition}
\begin{proof}
It follows from Lemma \ref{gcd lcm defined} that $P^{\rightlcm(g,g')}$ is
defined. It remains to show that
$k=\rightlcm(g,g')$ is $P$-reduced. Note that Lemma \ref{I-head preserved}
implies that
if $P\xrightarrow g$ is in the ribbon category and $g\preccurlyeq k$
then $g\preccurlyeq T_P(k)$. It follows that $g$ and $g'$ left-divide
 $T_P(k)$, thus $k=\rightlcm(g,g')$ left-divides $T_P(k)$ which proves
that $k$ is $P$-reduced.
\end{proof}
\begin{proposition}\label{Ribbon prefix}
For $P$ a standard parabolic submonoid, let $g\in M$ be $P$-reduced.
Then there is a unique maximal left-divisor of $g$ which is a $P$-ribbon.
If we denote by $R_P(g)$ this left-divisor, then
$R_P(g)\inv g$ is $P^{R_P(g)}$-reduced and there is equivalence between:
\begin{enumerate}
\item[(i)] $R_P(g)=1$.
\item[(ii)] Any atom which left-divides $\Delta_P g$ is in $P$.
\end{enumerate}
\end{proposition}
\begin{proof}
The  existence of  $R_P(g)$ is  a consequence  of the  fact that the ribbon
category  is  stable  by  right-lcms.  The  fact  that  $R_P(g)\inv  g$  is
$P^{R_P(g)}$-reduced is an immediate consequence of Lemma \ref{I-head preserved}.

We finally prove the equivalence of (i) and (ii) by observing that 
$R_P(g)\ne 1$ is equivalent to the existence
of an atom $s\notin P$ such that $v_{s,P}\preccurlyeq g$, which is in turn equivalent 
 to $\Delta_\Ps\preccurlyeq \Delta_P g$ which is equivalent to
$s$ left-dividing $\Delta_P g$.
\end{proof}
\begin{proposition} \label{positive conjugacy}
Let $P$ be standard parabolic submonoid of $M$, and let 
$b\in P$ and $g\in M$ be such that $g$ is  $P$-reduced and $b^g$
is defined. Assume
further that $\Delta_P\preccurlyeq b^i$ for some integer $i$. Then $g$ is
a $P$-ribbon.
\end{proposition}
\begin{proof}
The fact that $b^g$ is defined implies that $(b^i)^g$ is defined, thus replacing
$b$ by $b^i$ we may assume that $\Delta_P\preccurlyeq b$; let us thus write
$b=\Delta_P v$. 

The assumption that $b^g$ is defined can be written $g\preccurlyeq \Delta_Pvg$.
If $s$ is an atom left-dividing $g$ it follows, using Assumption
\ref{assumption}, that $\rightlcm(\Delta_P,s)=\Delta_\Ps\preccurlyeq
\Delta_Pvg$, whence $v_{s,P}\preccurlyeq vg$.
 By Lemma \ref{h<TP(g)} it follows that $v_{s,P}\preccurlyeq T_P(vg)=g$.

We conclude by Noetherian induction on $g$, since replacing simultaneously
$g$ by $v_{s,P}\inv g$, $P$ by $P^{v_{s,P}}$ and $b$ by
$b^{v_{s,P}}$ all the assumptions remain.
\end{proof}

We  conjecture that the assumption  in Proposition \ref{positive conjugacy}
that  $\Delta_P\preccurlyeq  b^i$  for  some  $i$  can  be  replaced by the
assumption that $P$ is the smallest standard parabolic submonoid containing
$b$.  For spherical Artin monoids, this is  a result of the beautiful paper
\cite{CGGW}. If we could prove that conjecture, we could extend all results
of \cite{CGGW} to our setting.

We now give a version of Proposition \ref{positive conjugacy} where
$g$ is not $P$-reduced.
\begin{proposition}\label{positive conjugate implies tail ribbon}
Let $P$ be a standard parabolic submonoid of $M$ and let $b\in P$, $g\in  M$ 
be such that $b^g\in M$.  
Assume that for some $i>0$ we have $\Delta_P\preccurlyeq b^i$.
Then $T_P(g)$ is a $P$-ribbon.
\end{proposition}
\begin{proof}
Let us prove first that $b^{H_P(g)}\in P$, that is $H_P(g)\preccurlyeq
bH_P(g)$.  Indeed, from $g\preccurlyeq bg$ we get $H_P(g)\preccurlyeq H_P(bg)
=bH_P(g)$.
Since $T_P(g)$ is $P$-reduced we can now
apply Proposition \ref{positive conjugacy}
with $b$ replaced by $b^{H_P(g)}$ and $g$
replaced by $T_P(g)$ and we get the result.
\end{proof}

\begin{lemma}\label{conjugating Delta_P conjugates P}
Let $P$ and $Q$ be standard parabolic submonoids of $M$ and $g\in M$, $k>0$ 
be such that $(\Delta_P^k)^g\in Q$;
then $T_P(g)$ is a $P$-ribbon and $P^{T_P(g)}\subset Q$.
If in addition $(\Delta_P^k)^g=\Delta_Q^k$ then $P^{T_P(g)}=Q$.
\end{lemma}
\begin{proof}
Proposition \ref{positive conjugate implies tail ribbon} applied with
$b=\Delta_P^k$ shows that $T_P(g)$ is a $P$-ribbon.
Replacing $\Delta_P^k$ by some power, we can assume that $k> l_\Delta(H_P(g))$, where
$l_\Delta$ is the Garside length, that is the number of
factors in a Garside normal form. We then have
$\Delta_P\preccurlyeq(\Delta_P^k)^{H_P(g)}\in P$. 
Now any atom of $P$ left-divides $(\Delta_P^k)^{H_P(g)}$ and is thus conjugate to an element of
$Q$ by $T_P(g)$, hence $P^{T_P(g)}\subset Q$.

If $(\Delta_P^k)^g=\Delta_Q^k$, let $Q'=P^{T_P(g)}$.
We have $\Delta_Q^k=(\Delta_P^k)^g=((\Delta_P^k)^{H_P(g)})^{T_P(g)}\subset
Q'$. Hence the divisors of $\Delta_Q$ in particular all the atoms of $Q$ are in $Q'$, so that $Q=Q'$.
\end{proof}
\section*{Garside Groups}
We now denote by $G$ the group of fractions of the Garside monoid $M$ (which
exists since $M$ is an Ore monoid, see \cite[2.32, V and 3.11, II]{livre}).

We call standard parabolic subgroup $G_P$ the group of fractions 
of a standard parabolic submonoid $P$.
\begin{definition} We say that $p\inv q$ is a left reduced
fraction  for $b\in G$ if $b=p\inv q$ with $p,q\in M$ and the
left-gcd of $p$ and $q$ is trivial.
\end{definition}
Symmetrically there are right reduced fractions $pq\inv$.
By ``reduced fraction'' we will mean left reduced fraction.

The reduced fraction for an element is unique; more precisely
if $p\inv q$ is reduced and $p\inv q=p^{\prime-1} q'$ there exists $d$
such that $p'=dp$ and $q'=dq$ (see \cite[3.11, II]{livre}).

By \cite[II, 3.18]{livre} $G_P$ is a subgroup of $G$ and $G_P\cap M=P$.
If an element of $G_P$ has $p\inv q$ as its reduced fraction in $G$, then
$p,q\in P$ since there is a reduced fraction $a\inv b$ in $P$ and 
$a\succcurlyeq p$ and $b\succcurlyeq q$. It follows that if $P$ and $Q$ are
standard parabolic submonoids and $P\subsetneq Q$ then $G_P\subsetneq G_Q$.
\begin{lemma}\label{Head in P}
Let $b\in G_P$ and $g\in M$ be such that $b^g\in M$. Then $b^{H_P(g)}\in P$.
\end{lemma}
\begin{proof}
Let $b'=b^g$, and write $g=H_P(g)T_P(g)$. Let
$b^{H_P(g)}=pq\inv$ where the right-hand side is a right reduced fraction in
$P$, that is  $p,q\in P$.
From the equality of the two fractions $T_P(g)b'T_P(g)\inv=pq\inv$
we deduce $q\preccurlyeq T_P(g)$ which implies $q=1$ since $T_P(g)$ is $P$-reduced.
\end{proof}
\begin{assumption}\label{pow delta central}
We assume that the automorphism induced by the conjugation by $\Delta$ is
of finite order, or equivalently that some power of $\Delta$ is central.
\end{assumption}
From now on we assume that Assumption \ref{pow delta central} holds. It
holds in spherical Artin monoids and in the dual braid monoids of well
generated complex reflection groups.
\begin{proposition}\label{conjugate into parabolic}
For two standard parabolic submonoids $P$ and $Q$ and $g\in G$ the three
following properties are equivalent
\begin{enumerate}
\item $G_P^g\subset G_Q$
\item There exists $k>0$ such that $(\Delta_P^k)^g\in G_Q$.
\item There exist $p\in G_P$ and a central power $\Delta^i$ of $\Delta$ such
that $P\xrightarrow{pg\Delta^i}$ is a ribbon such that $P^{pg}\subset Q$.
\end{enumerate}
\end{proposition}
\begin{proof}
Clearly (iii) implies (i) and (i) implies (ii). We prove that (ii)
implies (iii).
Let $u=(\Delta_P^k)^g\in G_Q$. 
Multiplying $g\inv$ by some central power $\Delta^i$ of $\Delta$ 
we get $h\in M$ such that $u^h=\Delta_P^k$. Then by
 the Lemma  \ref{Head in P} we have $u^{H_Q(h)}\in Q$. We thus have
$(\Delta_P^k)^{T_Q(h)\inv}\in Q$. Multiplying $T_Q(h)\inv$
by a central power $\Delta^j$
of $\Delta$ we get $m\in M$ such that $(\Delta_P^k)^m\in Q$. By
Lemma \ref{conjugating Delta_P conjugates P}  $P\xrightarrow{T_P(m)}$ is a
ribbon. Now $h\in g\inv \Delta^i$, thus $T_Q(h)\in
Qg\inv \Delta^i$, whence $m\in gG_Q\Delta^{j-i}$  and
$T_P(m)$ is in $P gG_Q \Delta^{j-i}=G_P g\Delta^{j-i}$, the last equality
since $g$ conjugates $G_Q$ to $G_P$.
\end{proof}
\begin{proposition}\label{conjugate parabolics}
Assume Assumption \ref{pow delta central}.
For two standard parabolic submonoids $P$, $Q$ we have equivalence between:
\begin{enumerate}
\item There is $g\in G$ such that $G_P^g\cap G_Q$ is of finite index in
$G_P^g$ and in $G_Q$.
\item There is $g\in G$ such that $G_P^g=G_Q$.
\item There exists an integer $k>0$ and $g\in G$ such that
 $(\Delta_P^k)^g=\Delta_Q^k$.
\item There exists a ribbon $P\xrightarrow g$ such that $P^g=Q$.
\end{enumerate}
\end{proposition}
\begin{proof}
It  is clear  that $g$  satisfying (iv)  satisfies (iii).  If $g$ satisfies
(iii),  we may assume that  $g\in M$ up to  multiplying $g$ by some central
power  of $\Delta$. Then Lemma \ref{conjugating Delta_P conjugates P} shows
that $T_P(g)$ satisfies (iv). Thus (iv) and (iii) are equivalent.

It  is  clear  that  $g$  satisfying  (iv)  satisfies  (ii),  and  that $g$
satisfying  (ii) satisfies (i). It is thus  enough to show that (i) implies
(iv).

So  we assume (i). Since $G_P^g\cap G_Q$ is of finite index in $G_P$, there
is some positive power $k$ such that $(\Delta_P^k)^g\in G_Q$. It follows by
Proposition \ref{conjugate  into parabolic}  that there  exists $g'\in  M$ of  the form
$g\Delta^i  q$ with  $\Delta^i$ a  central power  and $q\in  G_Q$ such that
$P\xrightarrow{g'}$ is a ribbon (and $P^{g'}\subset Q$).

The  element $g'$  satisfies clearly  the same  finite index assumptions as
$g$.  Thus there exists $k$ such that $G_P^{g'}\owns \Delta_Q^k$. Let $u\in
G_P$  be such  that $u^{g'}=\Delta_Q^k$.  By Lemma  \ref{Head in P} we have
$u^{H_P(g')}\in  P$, but  $H_P(g')=1$ thus  $P^{g'}\owns \Delta_Q^k$. This
implies  that  $P^{g'}=Q$  since  $P^{g'}$  is  a  standard parabolic
submonoid, thus
contains  all divisors  of its  element $\Delta_Q^k$,  that is all atoms of
$Q$.
\end{proof}
For a similar proposition for spherical Artin groups see \cite[Th\'eor\`eme 0.1]
{Godelle}.
\begin{corollary}\label{conjequiv}
It is equivalent for standard parabolic submonoids $P$ and $Q$ that $g\in G$ 
conjugates $G_P$ onto $G_Q$,
or that it conjugates some central power $\Delta_P^k$ to $\Delta_Q^k$.
If $g\in M$ is $P$-reduced it is equivalent that it conjugates $\Delta_P$ to
$\Delta_Q$ or that it conjugates $P$ to $Q$.

In particular if $G_P$ is conjugate to $G_Q$ and $\Delta_p^k$ is the smallest
power central in $P$ then $\Delta_Q^k$ is the smallest power central in $Q$.
\end{corollary}
\begin{proof}
We remark that in the proof (iii) $\Rightarrow$ (iv) $\Rightarrow$ (ii)
of Proposition \ref{conjugate parabolics}
the element obtained differs from $g$  by a central power of $\Delta$ and
an element of $P$, which
gives that if $g$ conjugates $\Delta_P^k$ to $\Delta_Q^k$ it conjugates
$P$ to $Q$.

In the proof of (ii) $\Rightarrow$ (iii)
of Proposition \ref{conjugate parabolics} the element obtained is 
in $G_P g\Delta^i G_Q$  for some central power $\Delta^i$. If $\Delta_P^k$
is central this element has the same effect as $g$ on it. This proves the
reverse implication.

The second sentence results from Lemma
\ref{conjugating Delta_P conjugates P}.
\end{proof}
We call parabolic subgroups the 
conjugates of the standard parabolic subgroups.  Note that the notion of
parabolic subgroup of $G$ depends on the Garside monoid $M$; in particular in
a spherical Artin monoid the parabolic subgroups for the ordinary Garside
structure are not the same as those for the dual Garside structure.
\begin{definition}
If  $K=G_P^g$  is  a  parabolic subgroup, we
denote  by $z_K$ the element $(\Delta_P^k)^g$ where $\Delta_P^k$ is the
smallest central power of $\Delta_P$.
\end{definition}
The notation $z_K$ makes sense thanks to the following proposition:
\begin{proposition}  Let $K,K'$ be parabolic subgroups of $G$.
\begin{itemize} 
\item $z_K$ depends only on $K$.
\item It is equivalent that $g\in G$ conjugates $K$ to $K'$ or $z_K$ to
$z_{K'}$.
\end{itemize} 
\end{proposition}
\begin{proof}
For the first item, let $K=G^g_P$ and $K=G^{g'}_{P'}$ two ways of conjugating
$K$ to a standard parabolic subgroup. This defines two candidates for $z_K$
which are $(\Delta_P^k)^g$ and $(\Delta_{P'}^k)^{g'}$. But these two elements
are equal by Corollary \ref{conjequiv} since $g'g\inv$ conjugates $P'$ to
$P$, thus $\Delta_{P'}^k$ to $\Delta_P^k$.

For the second item, let $K=G^g_P$ and $K'=G^{g'}_{P'}$ be ways of conjugating
$K,K'$ to standard parabolic subgroups. Then $z_K=(\Delta_P^k)^g$ and 
$z_{K'}=(\Delta_{P'}^k)^{g'}$. Now $K^x=K'$ for some $x\in G$
is equivalent to $G_P^{gxg^{\prime-1}}=
G_{P'}$ and we conclude again by Corollary \ref{conjequiv}.
\end{proof}
\begin{proposition}\label{standardizer}
Let $K=G_P^b$ be a parabolic subgroup of $G$, where $P$ is a standard
parabolic submonoid and $b\in M$. 
Define  $b'$ by
$H_P(b)R_P(T_P(b))b'=b$ and $Q$ by $Q=P^{R_P(T_P(b))}$.
Then $b^{\prime -1}\Delta_Q^k b'$ is the reduced fraction of $z_K$, where
$\Delta_Q^k$ is the smallest central power of $\Delta_Q$ in $Q$.
\end{proposition}
\begin{proof}
We first remark that by definition we have $(\Delta_P^k)^b=z_K$. We may clearly
replace $b$ in this equality by $T_P(b)$. Let $c=R_P(T_P(b))$; we have
$(\Delta_P^k)^c=\Delta_Q^k$. We thus get $z_K=b^{\prime-1}\Delta_Q^k b'$.
We claim this is a reduced fraction. Indeed by construction $R_Q(b')=1$ thus
by Proposition \ref{Ribbon prefix}(ii) any atom left-dividing
$\Delta_Q b'$ is in $Q$,
thus the same is true for $\Delta_Q^k b'$ by induction, using that
for $k>1$ one has $H(\Delta_Q^k b')=H(\Delta_Q H(\Delta_Q^{k-1}b'))$.
Since $b'$ is $Q$-reduced, the fraction is reduced.
\end{proof}
Note that $b'$ above is  minimal such that $\lexp{b'}K$ is standard, that is
any $u\in M$ such that $\lexp uK$ is 
standard is a left multiple of $b'$, and $G_Q$ is a ``canonical''
standard parabolic subgroup conjugate to $K$. For spherical Artin groups
the proposition is \cite[Theorem 3]{Cumplido} and the
element $b'$ is called a minimal standardizer.

An immediate consequence of Proposition \ref{standardizer} if that it is
equivalent that $z_K\in M$ or that $K$ is standard.


\begin{thebibliography}{x}
\bibitem{Bessis}D.~Bessis, Finite complex reflection arrangements are 
$K(\pi,1)$, {\sl Ann.~of Math.~\bf 181} (2015) 809--904 

\bibitem{Cumplido}M.~Cumplido, On the  minimal positive standardizer of a 
parabolic subgroup of an Artin-Tits group, 
{\sl J. Algebraic Combinatorics \bf 49} (2019) 337--359

\bibitem{CGGW}M.~Cumplido, V.~Gebhardt, J.~Gonzalez-Meneses and B.~Wiest,
On parabolic subgroups of Artin-Tits groups of spherical type,
{\sl Adv. Math. \bf 352} (2019) 572--610

\bibitem{Godelle}E.~Godelle, Normalisateur et groupe d'Artin de type
sph\'erique, {\sl J. Algebra \bf 269} (2003) 263--274

\bibitem{livre}P.~Dehornoy, F.~Digne, E.~Godelle, D.~Krammer, J.~Michel,
Foundations of Garside theory, {\sl EMS Tracts in Math.~\bf 22} (2015)

\bibitem{DMM} F.~Digne, I.~Marin and J.~Michel,
The center of pure complex braid groups
{\sl J.~Algebra \bf 347} (2011) 206--213 
\end{thebibliography}
\end{document}